\input{amssymb.sty}
\documentclass{amsart}
\theoremstyle{plain}
\newtheorem{thm}{Theorem}[section]
\newtheorem*{thm*}{Theorem}
\newtheorem{cor}[thm]{Corollary}
\newtheorem*{cor*}{Corollary}
\newtheorem*{conj*}{Conjecture}
\newtheorem*{lemma*}{Lemma}
\newtheorem{lemma}[thm]{Lemma}
\newtheorem*{prop*}{Proposition}
\newtheorem{prop}[thm]{Proposition}

\theoremstyle{definition}
\newtheorem{rems}[thm]{Remarks}
\newtheorem*{defn*}{Definition}
\newtheorem*{rems*}{Remarks}
\newtheorem*{proof*}{Proof}
\newtheorem{prel*}{Preliminaries}
\newtheorem{examples*}{Examples}

\newcommand{\C}{\mathbb{C}}

\newcommand{\D}{{\mathbb D}}

\newcommand{\Ch}{\operatorname{Ch}}

\newcommand{\Index}{\operatorname{index}}

\newcommand{\tr}{\operatorname{tr}}

\newcommand{\R}{\operatorname{\mathbb R}}
\newcommand{\Z}{\operatorname{\mathbb Z}}

\newcommand{\nc}{\newcommand}
\nc{\nt}{\newtheorem}
\nc{\gf}[2]{\genfrac{}{}{0pt}{}{#1}{#2}}
\nc{\mb}[1]{{\mbox{$ #1 $}}}
\nc{\real}{{\mathbb R}}
\nc{\comp}{{\mathbb C}}
\nc{\ints}{{\mathbb Z}}
\nc{\Ltoo}{\mb{L^2({\mathbf H})}}
\nc{\rtoo}{\mb{{\mathbf R}^2}}
\nc{\slr}{{\mathbf {SL}}(2,\real)}
\nc{\slz}{{\mathbf {SL}}(2,\ints)}
\nc{\su}{{\mathbf {SU}}(1,1)}
\nc{\so}{{\mathbf {SO}}}
\nc{\hyp}{{\mathbb H}}
\nc{\disc}{{\mathbf D}}
\nc{\torus}{{\mathbb T}}

\nc{\ca}{{\mathcal A}}
\nc{\cag}{{{\mathcal A}^\Gamma}}
\nc{\cg}{{\mathcal G}}
\nc{\chh}{{\mathcal H}}
\nc{\ck}{{\mathcal B}}
\nc{\cl}{{\mathcal L}}
\nc{\cm}{{\mathcal M}}
\nc{\cs}{{\mathcal S}}
\nc{\cz}{{\mathcal Z}}
\nc{\sind}{\sigma{\rm -ind}}

\begin{document}

\title[QUANTUM HALL EFFECT ON THE HYPERBOLIC PLANE..] {QUANTUM HALL EFFECT ON
THE HYPERBOLIC PLANE IN THE PRESENCE OF DISORDER} \author{A. Carey}
\address{Department of Mathematics, University of Adelaide, Adelaide 5005,
Australia}
\email{acarey@maths.adelaide.edu.au}
\author{K. Hannabuss}
\address{Department of Mathematics, University of Oxford, England.}
\email{khannabu@maths.adelaide.edu.au}
\author{V. Mathai}
\address{Department of Mathematics, University of Adelaide, Adelaide 5005,
Australia}
\email{vmathai@maths.adelaide.edu.au}


\begin{abstract} We study both the continuous model and the discrete model of
the quantum Hall effect (QHE) on the hyperbolic plane in the presence of
disorder, extending the results of an earlier paper \cite{CHMM}. Here we model
impurities, that is we consider the effect of a {\em random or almost periodic
potential} as opposed to just periodic potentials. The Hall conductance is
identified as a geometric invariant associated to an algebra of observables,
which has plateaus at {\em gaps in extended states} of the Hamiltonian. We use
the Fredholm modules defined in \cite{CHMM} to prove the integrality of the
Hall
conductance in this case. We also prove that there are always only a finite
number of gaps in extended states of any random discrete Hamiltonian.
\end{abstract}

\maketitle

\section*{Introduction}

In \cite{CHMM}, we studied continuous and discrete models for the integer
quantum Hall effect on the hyperbolic plane, generalising the study of the
Euclidean model due to Bellissard \cite{Bel+E+S} and Xia \cite{Xia}. Our model
involves electrons moving in a two dimensional conducting material with
hyperbolic geometry. One should think of hyperbolic space and hence the sample
as an embedded hyperboloid in pseudo-Euclidean 3-space. The crystal lattice of
the conductor is now modelled by the orbit of a freely acting discrete group
$\Gamma$. For reasons of convenience we take this to be the fundamental
group of
a Riemann surface (though aspects of our analysis work more generally). The
magnetic field remains orthogonal to the two dimensional conductor but we did
not attempt to model impurities and so we allowed only $\Gamma$ invariant
potential terms. No assumption was made about the rationality of the imposed
magnetic flux while the integrality of the Hall conductance follows by showing
that it is given by the index of a Fredholm operator (we also exhibited the
conductance as a topological index). Thus the models in \cite{CHMM} were
noncommutative geometry models for the quantum Hall effect on hyperbolic
space.

In this letter, we extend the results in \cite{CHMM} in two different
directions. First, we now model impurities, in other words, we generalise to
almost periodic or random potentials. Second, we extend the definition of Hall
conductance to projections to gaps in extended states of the Hamiltonian. These
include the projections in \cite{CHMM} but are potentially much richer. We also
show that the number of gaps in extended states for the discrete model is
always
finite. Some of the results of this paper are based on unpublished work
\cite{Ma} that was announced in \cite{CHMM}.

We begin by reviewing the construction of the Hamiltonian in the continuous
model for the quantum Hall effect on hyperbolic space. First, we take as our
principal model of hyperbolic space, the upper half-plane $\hyp$ in $\comp$
equipped with its usual Poincar\'e metric $(dx^2+dy^2)/y^2$, and symplectic
area
form $\omega_\hyp = dx\wedge dy/y^2$. The group $\slr$ acts transitively on
$\hyp$ by M\"obius transformations $$ x+iy = \zeta \mapsto g\zeta =
\frac{a\zeta+b}{c\zeta+d},\quad\mbox{for } g=\left(
\begin{array}{cc} a & b\\ c & d
\end{array}\right).$$ Any Riemann surface of genus $g$ greater than 1 can be
realised as the quotient of
$\hyp$ by the action of its fundamental group realised as a subgroup
$\Gamma$ of
$\slz \subset \slr$.

Pick a 1-form $\eta$ such that $d\eta = \theta\omega_\hyp$, for some fixed
$\theta \in \R$. As in geometric quantisation we may regard $\eta$ as
defining a
connection $\nabla = d-i\eta$ on a line bundle $\cl$ over $\hyp$, whose
curvature is $\theta\omega_\hyp$. Physically we can think of $\eta$ as the
electromagnetic vector potential for a uniform magnetic field of strength
$\theta$ normal to $\hyp$. Using the Riemannian metric the Hamiltonian of an
electron in this field is given in suitable units by $$H = H_\eta = \frac
12\nabla^*\nabla = \frac 12(d-i\eta)^*(d-i\eta).$$ Comtet [Comtet] has shown
that $H$ differs from a multiple of the Casimir element for $\slr$, $\frac
18{\bf J}.{\bf J}$, %+\frac 14B^2$ by a constant, where $J_1$, $J_2$ and $J_3$
denote a certain representation of generators of the Lie algebra $sl(2,\real)$,
satisfying $$[J_1,J_2] = -iJ_3, \qquad [J_2,J_3] = iJ_1, \qquad [J_3,J_1] =
iJ_2,$$ so that ${\bf J}.{\bf J} = J_1^2+J_2^2-J_3^2$ is the quadratic Casimir
element. This shows very clearly the underlying $\slr$-invariance of the
theory.
Comtet has computed the spectrum of the unperturbed Hamiltonian $H_\eta$, for
$\eta = -\theta dx/y$, to be the union of finitely many eigenvalues
$\{(2k+1)\theta -k(k+1):k=0,1,2\ldots < \theta-\frac 12\}$, and the continuous
spectrum $[\frac 14 + \theta^2, \infty)$. Any $\eta$ is cohomologous to
$-\theta
dx/y$ (since they both have $\omega_\hyp$ as differential) and forms differing
by an exact form $d\phi$ give equivalent models: in fact, multiplying the wave
functions by $\exp(i\phi)$ shows that the models for $\eta$ and $-\theta dx/y$
are unitarily equivalent. This equivalence also intertwines the
$\Gamma$-actions
so that the spectral densities for the two models also coincide.

A more realistic model Hamiltonian is obtained by adding a potential $V$. In
\cite{CHMM}, we took
$V$ to be invariant under $\Gamma$. Here, in section 1, we model impurities by
allowing any smooth random potential function $V$ on $\mathbb H$. We consider
two general notions of  random potential whereas in the literature the
$\Gamma$-action on the disorder space is required to admit an ergodic invariant
measure. The class of random potentials we consider contains all the periodic
potentials, but it is much larger. In particular, we will show that any smooth
bounded potential $V$ is a random potential for some disorder space $\Omega$.
The perturbed Hamiltonian $H_{\eta,V} = H_\eta +V$ has unknown spectrum for
general smooth and bounded $V$. However we are able to deduce some qualitative
aspects of the spectrum of these Hamitonians by using a reduction (via Morita
equivalence) to a simpler case: that of a discrete model. In Section 2 we
review
the construction of the Hamiltonian in this discrete model for the quantum Hall
effect on hyperbolic space again taking, as our principal model, the hyperbolic
plane. The unperturbed Hamitonian
$H_\sigma$ is then the random walk operator on the Cayley graph of $\Gamma$ in
the projective $(\Gamma, \sigma)$-representation that is defined by the
magnetic
field. To model impurities this Hamiltonian is modified by the addition of a
perturbation $T$ which is almost periodic or random. Again, it is the perturbed
Hamiltonian $H_{\sigma, T} = H_\sigma +T$ whose spectral properties are of
interest and about which we can obtain qualitative information.

In section 3, we derive the hyperbolic Connes-Kubo formula for the Hall
conductance for random Hamiltonians in the discrete model, following \cite{MM}.
We show that even in this case, the hyperbolic Connes-Kubo formula is
cohomologous to another cyclic 2-cocycle, from which we can easily deduce
integrality by our previous results in \cite{CHMM}. We then deduce our first
version of our theorem concerning the quantum Hall effect on the hyperbolic
plane in the presence of disorder, showing that the plateaus of the Hall
conductance occur at the gaps in the spectrum of the random Hamiltonian, where
the Hall conductance takes on integral values in $2(g-1)\Z$.

In section 4, we study the domain of definition of these cyclic 2-cocycles. The
gaps in extended states of a random Hamiltonian correspond precisely to the
plateaus of the Hall conductance in the quantum Hall effect. Our next goal
is to
prove that there are only a {\em finite} number of gaps in extended states of
the random Hamiltonian $H_{\sigma, T}$, exactly corresponding to what is
expected to be physically observed. In contrast, it is sometimes possible to
have an infinite number of gaps in the spectrum of the random Hamiltonian
$H_{\sigma, T}$, cf.
\cite{CHMM}, \cite{CEY}. Our main result is found in the last section, where we
show that the plateaus of the Hall conductance occur at the gaps in extended
states of the random Hamiltonian whose disorder space has an ergodic
$\Gamma$-action, with the Hall conductance taking on integral values in
$2(g-1)\Z$.

\section{Continuous model}

\subsection{The geometry of the hyperbolic plane}

The upper half-plane can be mapped by the Cayley transform $z =
(\zeta-i)/(\zeta+i)$ to the unit disc $\D$ equipped with the metric
$|dz|^2/(1-|z|^2)^2$ and symplectic form $dz\,d\overline{z}/2i(1-|z|^2)^2$, on
which $\su$ acts, and some calculations are more easily done in that
setting. In
order to preserve flexibility we shall work more abstractly with a Lie
group $G$
acting transitively on a space $X \sim G/K$. Although we shall ultimately be
interested in the case of $G = \slr$ or $\su$, and $K$ the maximal compact
subgroup which stabilises $\zeta = i$ or $z=0$ so that $X=\mathbb H$ or $X=
\mathbb D$, those details will play little role in many of our calculations,
though we shall need to assume that $X$ has a $G$-invariant Riemannian metric
and symplectic form $\omega_X$. We shall denote by $\Gamma$ a discrete subgroup
of $G$ which acts freely on $X$ and hence intersects $K$ trivially.

We shall assume that $\cl$ is a hermitian line bundle over $X$, with a
connection, $\nabla$, or equivalently, for each pair of points $w$ and $z$ in
$X$, we denote by $\tau(z,w)$ the parallel transport operator along the
geodesic
from $\cl_w$ to $\cl_z$. In $\hyp$ with the line bundle trivialised and $\eta =
\theta dx/y$ one can calculate explicitly that $$\tau(z,w) =
\exp\left(i\int_w^z\eta\right) = [(z-\overline{w})/(w-\overline{z})]^\theta.$$
For general $\eta$ we have $\eta - \theta dx/y = d\phi$ and $$\tau(z,w) =
\exp(i\int_w^z\eta) =
[(z-\overline{w})/(w-\overline{z})]^\theta\exp(i(\phi(z)-\phi(w))).$$ Parallel
transport round a geodesic triangle with vertices $z$, $w$, $v$, gives rise
to a
holonomy factor:
$$\varpi(v,w,z) = \tau(v,z)^{-1}\tau(v,w)\tau(w,z),$$ and this is clearly the
same for any other choice of $\eta$, so we may as well work in the general
case.

\begin{lemma} The holonomy can be written as
$\varpi(v,w,z) = \exp\left(i\theta\int_\Delta \omega_\hyp\right)$, where
$\Delta$ denotes the geodesic triangle with vertices $z$, $w$ and $v$. The
holonomy is invariant under the action of $G$, that is $\varpi(v,w,z) =
\varpi(gv,gw,gz)$, and under cyclic permutations of its arguments.
Transposition
of any two vertices inverts $\varpi$. For any four points $u$ ,$v$, $w$, $z$ in
$X$ one has $$\varpi(u,v,w)\varpi(u,w,z) = \varpi(u,v,z)\varpi(v,w,z).$$
\end{lemma} \begin{proof} By definition, for a suitable trivialisation of $\cl$
one has $$\varpi(v,w,z) = \exp\left(i\int_{\partial\Delta}\eta\right)$$ and the
first part follows by applying Stokes' Theorem after noting that the result is
independent of the trivialisation. The invariance under $G$ follows from the
invariance of the symplectic form, and the results of permutations follow from
the properties of the integral, as does the final identity. \end{proof}

\subsection{Algebra of observables and random or almost periodic
potentials} The
algebra of physical observables that we consider in the continuous model should
include the operators $f(H_{\eta, V})$ for any bounded continuous function $f$
on $\R$ and for any smooth random potential function
$V$ on $\mathbb H$. The class of random potentials contains all the periodic
potentials, but it is much larger. In particular, we will show that any smooth
bounded potential $V$ is a random potential for some disorder space $\Omega$.

We will show that the twisted $C^*$-algebra of the groupoid
$\mathcal G = \Gamma\backslash (X \times X \times \Omega)$, twisted by
$\varpi$,
is large enough to contain all the physical observables. This algebra also
turns
out to be the twisted $C^*$-algebra  of the foliation $\Omega_\Gamma$. We will
also show that this $C^*$-algebra is strongly Morita equivalent to the cross
product $C^*$-algebra $C(\Omega) \rtimes_\sigma \Gamma$, where $\sigma$ is a
multiplier on $\Gamma$ which is determined by $\varpi$.

The assumptions that we impose on the {\em disorder space} $\Omega$ are:
\begin{itemize}
\item[(1)] $\Omega$ is compact and admits a Borel probability measure
$\Lambda$;
\item[(2)] there is a continuous action of $\Gamma$ on $\Omega$ with a dense
orbit. \end{itemize}

The geometrical data described in the last subsection enables us to easily
describe the first of the two $C^*$algebras which appear in the theory. This
twisted algebra of kernels, which was introduced by Connes \cite{Co2} is the
$C^*$-algebra $\ck$ generated by compactly supported smooth functions on
$X\times X \times \Omega$ with the multiplication $$k_1*k_2(z,w, r) = \int_X
k_1(z,v, r)k_2(v,w, r)\varpi(z,w,v)\,dv,$$ (where $dv$ denotes the
$G$-invariant
measure defined by the metric) and $k^*(z,w, r)= \overline{k(w,z, r)}$.
There is
an obvious trace on $\ck$ given by $\tau_\ck (k) = \int_{X\times \Omega} k(z,z,
r)\,dz d\Lambda(r)$. Observe that $X\times X\times\Omega$ is a groupoid with
space of units $X\times \Omega$ and with source and range maps $s((z,w,r)) =
(w,r)$ and $r((z,w,r')) = (z,r')$. Then the algebra of twisted kernels is the
extension of the $C^*$-algebra of the groupoid $X\times X\times\Omega$ defined
by the cocycle $((v,w, r),(w,z, r)) \mapsto \varpi(v,w,z)$, \cite{Ren1}.

\begin{lemma} The algebra $\ck$ has a representation $\pi$ on the space of
$L^2$
sections of $\cl \to X\times \Omega$ defined by
$$(\pi(k)\psi)(z, r) = \int_X k(z,w, r)\tau(z,w)\psi(w, r)\,dw.$$ \end{lemma}
\begin{proof} The parallel transport $\tau(z,w)$ ensures that the integral
is in
the appropriate fibre, and the fact that it is a representation follows from a
calculation using the definition of the holonomy. \end{proof}

We now pick out a
$\Gamma$-invariant subalgebra $\ck^\Gamma$ of $\ck$. This condition reduces
simply to the requirement that the kernel satisfies
$k(\gamma^{-1}z,\gamma^{-1}w, \gamma^{-1}r) = k(z,w, r)$ for all $\gamma\in
\Gamma$. As before, observe that $\Gamma\backslash(X\times X\times\Omega)$ is a
groupoid whose elements are $\Gamma$ orbits $(x,y, v)_\Gamma = \{(\gamma
x,\gamma y, \gamma v): \gamma \in \Gamma\}$, with source and range maps
$s((x,y, v)_\Gamma) = ( y, v)_\Gamma$ and $r((x,y, v)_\Gamma) = ( x,
v)_\Gamma$.
The space of units is $\Omega_\Gamma = \Gamma\backslash(X\times \Omega)$, whose
elements are $\Gamma$ orbits $(x, v)_\Gamma = \{(\gamma x, \gamma v):
\gamma \in
\Gamma\}$. Then the algebra of invariant twisted kernels $\ck^\Gamma$ is the
extension of the $C^*$-algebra of the groupoid $\Gamma\backslash(X\times
X\times\Omega)$ defined by the cocycle $((v,w, r),(w,z, r)) \mapsto
\varpi(v,w,z)$, \cite{Ren1}. With our assumptions on the disorder space
$\Omega$, there is in general {\em no} trace on the algebra $\ck^\Gamma$, and
there may not even be a weight on this algebra in general. However, we mention
that under the additional assumption that the measure $\Lambda$ on $\Omega$ is
$\Gamma$-invariant, the natural trace $\tau_{\ck^\Gamma}$ for this algebra is
given by the same formula as before except that the integration is now over
$\Omega_\Gamma = \Gamma \backslash (X\times \Omega)$ rather than $X\times
\Omega$, where we have identified $\Omega_\Gamma$ with a fundamental domain: $
\tau_{\ck^\Gamma} (T)=\int_{\Omega_\Gamma} T(z,z, r)dzd\Lambda(r). $ We also
mention that under the additional assumption that the measure $\Lambda$ on
$\Omega$ is quasi-$\Gamma$-invariant, the natural tracial weight
$\tau_{\ck^\Gamma}$ for this algebra is given by $
\tau_{\ck^\Gamma} (T)=\int_{X\times\Omega} f(z,r)^2T(z,z, r)dzd\Lambda(r), $
where $f\in C_c(X\times\Omega)$ is such that $\sum_{\gamma\in \Gamma}
(\gamma^*f)^2 = 1$.

We now recall a notion due to Connes \cite{Co2}.

\begin{defn*} A {\em random} or {\em almost periodic} potential on $X$ is a
continuous family of smooth functions on the disorder space, $ \Omega \ni r
\mapsto V_r \in C^\infty(X)$ where the following equivariance is imposed: $$
V_{\gamma r} = \gamma^* V_r \qquad \forall \gamma \in \Gamma, \forall r\in
\Omega.
$$
\end{defn*}

\begin{rems} If $V$ is a $\Gamma$-invariant potential on $X$, then it is
clearly
random for any disorder space. More generally, if $V$ is a arbitrary smooth
function on $X$ such that the set $\left\{ \gamma^* V : \gamma \in
\Gamma\right\}$ has compact closure in the strong operator topology in
$B(L^2(X))$, then $V$ is a random potential. \end{rems}

We have asserted informally that the Hamiltonian can be accommodated within the
algebra $\ck^\Gamma$ and we shall now provide the proof. Fix a base point
$u\in\D$ and introduce: $$\sigma(x,y) = \varpi(u,xu,xyu)$$
$$\phi(z,\gamma) = \varpi(u, \gamma^{-1}u,\gamma^{-1}z)\tau(u,z)^{-1}
\tau(u,\gamma^{-1}z).$$ Then $\sigma$ is the group 2-cocycle in the projective
action of $\su$ on $L^2(\D)$ defined by:
$$U(\gamma)\psi(z) = \phi(z,\gamma)\psi(\gamma^{-1}z)$$ where
$\psi\in L^2(\D), \gamma\in\su$. Note that $U$ is constructed so that the
$\Gamma$-invariant algebra
$\pi(\ck^\Gamma)$ is the intersection of $\pi(\ck)$ with the commutant of $U$.
Recall that the unperturbed Hamiltonian $H = H_\eta$ commutes with the
projective representation $U$ (cf. Lemma 4.9, \cite{CHMM}). So we see that $H$
is affiliated to the von Neumann algebra generated by the representation $\pi$
of $\ck^\Gamma$ (cf. Corollary 4.2 \cite{CHMM}).

A random potential $V$ can be viewed as defining an equivariant family of
Hamiltonians $\Omega \ni r \mapsto H_{\eta, V_r} = H + V_r \in {\rm
Oper}(L^2(X)$ where ${\rm Oper}(L^2(X))$ denotes closed operators on $L^2(X)$.

Br\"uning and Sunada have proved an estimate on the Schwartz kernel of the heat
operator for any elliptic operator, and in particular for $\exp(-t H_{\eta,
V_r})$ for $t > 0$, which implies that it is $L^1$ in each variable separately.
Since this kernel is $\Gamma$-equivariant it follows (in exactly the same
fashion as Lemma 4 of \cite{BrSu}) that this estimate implies that $\exp(-t
H_{\eta, V_r})$ is actually in the algebra $\ck^\Gamma$.

\begin{lemma} One has
$f(H_{\eta, V}) \in \ck^\Gamma$ for any bounded continuous function $f$ on $\R$
and for any random potential $V$ on $X$. In particular, the spectral
projections
of $H_{\eta, V}$ corresponding to gaps in the spectrum lie in $\ck^\Gamma$.
\end{lemma} \begin{proof} If $\mu$ lies in a gap of the spectrum of $H_{\eta,
V}$ then let $f$ be a continuous approximate step function which is identically
one on the part of the spectrum of $H$ contained in $[0,\mu]$ and zero on the
part contained in $[\mu,\infty)$. Define $g(x)=f\circ\ln(1/x)$ for $x\in
(0,1]$.
Then $g$ is a bounded continuous function which, when applied to $e^{-tH}$
gives
the spectral projection corresponding to the interval $[0,\mu]$. More
generally,
any bounded continuous function $f$ on $\R$ can be uniformly approximated by
step functions, which proves the lemma. \end{proof}

The following is adapted from \cite{Bel+E+S},\cite{Nak+Bel}. Note that our
assumptions on the disorder space $\Omega$ are weaker than theirs.

\begin{thm} Let $V$ be a smooth bounded function on $X$. Then $V$ is a random
potential for some disorder space $\Omega$ and therefore $f(H_{\eta, V}) \in
\ck^\Gamma$ for any bounded continuous function $f$ on $\R$. \end{thm}

\begin{proof} Firstly, observe that the set $\left\{ \gamma^* V : \gamma \in
\Gamma\right\}$ is contained in the ball $$\left\{ V_1 \in L^\infty(X) :
||V_1||_\infty \le ||V||_\infty \right\}.$$ Also, any ball in $ L^\infty(X)$ is
compact in the weak topology, so the set $\left\{ \gamma^* V : \gamma \in
\Gamma\right\}$ also has compact closure in the weak topology in $
L^\infty(X)$.
Equivalently, the set $\left\{ \gamma^* V : \gamma \in \Gamma\right\}$ has
compact closure in the strong operator topology in $B(L^2(X))$, proving
that $V$
is a random potential. The proof is completed by now applying the previous
Lemma. \end{proof}

\noindent{\em Example} Let the Iwasawa decomposition of $\su$ be written $KAN$
then
$PSL(2,\Z)$ acts on $\D=\su/K$ by M\"obius transformations so that
$\Gamma\subset PSL(2,\Z)$ also acts. Let
$g_{\lambda,w}(z)=\lambda\frac{1-|z|^2}{|w-z|^2}$ where $\lambda\in \R^+\cong
A$, $w\in U(1)\cong K$ and $z\in\D$. Now let $\gamma=\left(
\begin{array}{cc}
\alpha & \beta\\
\bar\beta & \bar\alpha
\end{array}\right)$ and we calculate
$$ U(\gamma)g_{\lambda,w}U(\gamma^{-1})=g_{\lambda_{\gamma,w}\lambda,\gamma
w} $$
where $\lambda_{\gamma,w}= |\bar\beta w +\bar\alpha|^{-2}$. The stabiliser of
$g_{1,1}$ is
$\{\pm\left(
\begin{array}{cc} 1-in & in\\ -in & 1+in
\end{array}\right):\ n\in\R\}.$ This group is $MN$ where $MAN$ is the maximal
parabolic subgroup. Thus we have the usual action of $\su$ on $\su/MN$ and
hence
{\em a fortiori} a $\Gamma$-action which is known to be ergodic, cf.
\cite{Zim}.
Note that, regarding $\{e^{-g_{\lambda,w}}\}$ as a set of bounded
multiplication
operators on $L^2(\D)$, the strong closure of
$\{U(\gamma)e^{-g_{\lambda,w}}U(\gamma^{-1})\ |\ \lambda\in\R, w\in U(1)\}$ is
homeomorphic to $S^2$. (This is because taking the strong closure adds the zero
and identity operator to the set.) Thus in this example the disorder space is
$S^2$ which admits a dense orbit and a quasi-invariant ergodic probability
measure.

\subsection{Morita equivalence} We will prove here that the algebra
$\ck^\Gamma$
is Morita equivalent to another more tractable algebra. We shall do this by
using the groupoid equivalence arguments of \cite{M+R+W}, or rather the twisted
version, \cite{Ren2}, \cite{Ren3}. We have already noted that $\ck$ is the
$C^*$-algebra of an extension of the groupoid $X\times X\times \Omega$ by a
cocycle defined by $\varpi$, and $\Gamma$ invariance of $\varpi$ means that
$\ck^\Gamma$ is likewise the $C^*$-algebra of an extension of
$\Gamma\backslash(X\times X\times \Omega)$ by $\varpi$, where
$\Gamma\backslash(X\times X\times \Omega)$ denotes the groupoid obtained by
factoring out the diagonal action of $\Gamma$. More precisely, the groupoid
elements are $\Gamma$ orbits $(x,y, v)_\Gamma = \{(\gamma x,\gamma y,
\gamma v):
\gamma \in \Gamma\}$, with source and range maps
$s((x,y, v)_\Gamma) = ( y, v)_\Gamma$ and $r((x,y, v)_\Gamma) = ( x,
v)_\Gamma$. Therefore
$(x_1,y_1, v_1)_\Gamma$ and $(x_2,y_2, v_2)_\Gamma$ are composable if and only
if $y_1 = \gamma x_2$ and $v_1 = \gamma v_2$ for some $\gamma \in \Gamma$, and
then the composition is $(x_1,\gamma y_2, \gamma v_2)_\Gamma$. We also note
that
$\Omega \times \Gamma$ is a groupoid. The source and range maps are $s((v,
\gamma)) = \gamma v$ and $r((v, \gamma)) = v$. Therefore the elements
$(v_1, \gamma_1)$ and $(v_2, \gamma_2)$ are composable if and only if $v_1 =
\gamma_2 v_2$, and the composition is $(\gamma_2^{-1} v_1, \gamma_1 \gamma_2)$.

\begin{thm} The algebra $\ck^\Gamma$ is Morita equivalent to the twisted cross
product algebra
$C(\Omega) \rtimes_{\bar\sigma} \Gamma$. \end{thm}

\begin{proof} This result will follow immediately from \cite{Ren2} Corollaire
5.4 (cf \cite{M+R+W} Theorem 2.8) once we have established the groupoid
equivalence in the following Lemma.
\end{proof}

\begin{lemma} The line bundle $\cl$ over $X\times \Omega$ provides an
equivalence (in the sense of \cite{Ren2} Definition 5.3) between the groupoid
extensions $(\Gamma\backslash(X\times X\times \Omega))^\varpi$ of
$\Gamma\backslash(X\times X\times \Omega)$ defined by $\varpi$ and
$(\Omega \times\Gamma)^\sigma$ of $\Omega \times\Gamma$ defined by
$\overline{\sigma}$.
\end{lemma}

\begin{proof} Both extensions are by $\torus$. We write the elements of
$(\Gamma\backslash(X\times X\times \Omega))^\varpi$ as quadruples
$(x,y,v,t) \in
X\times X\times \Omega\times\torus$ with the first three elements
representing a
diagonal $\Gamma$ orbit. Elements are composable if their first three
components
are composable, and, when $y_1=\gamma x_2$ and $v_1 = \gamma v_2$, one has
$$(x_1,y_1,v_1, t_1)(x_2,y_2,v_2, t_2) = (x_1,\gamma y_2, \gamma v_2,
t_1t_2\varpi(x_1,y_1,\gamma y_2)).$$

The line bundle can be trivialised and written as $X\times \Omega \times\comp$,
since it is a pullback of a line bundle from $X$. We let
$(\Gamma\backslash(X\times X\times \Omega))^\varpi$ act on the left of the line
bundle by defining $(x,y,v,t)$ to act on $(z,v',u)$ if $z = \gamma y$ and
$\gamma v = v'$ for some
$\gamma \in \Gamma$, and then the result of the action is $(\gamma x,\gamma v,
\tau(\gamma x,z)tu)$. (One may check that this gives an action using the
relationship between parallel transport and holonomy and the
$\Gamma$-invariance
of $\varpi$.)

The twisted groupoid $(\Omega\times \Gamma)^\sigma$ has underlying set
$\Omega\times \Gamma\times \torus$. The elements $(v_1, \gamma_1, t_1)$ and
$(v_2, \gamma_2, t_2)$ are composable if and only if $v_1 = \gamma_2 v_2$, and
the composition is, $$(\gamma_2^{-1} v_1,
\gamma_1\gamma_2,\sigma(\gamma_1,\gamma_2)t_1 t_2).$$ Next, $(\Omega\times
\Gamma)^\sigma$ acts on the right of $\cl$ whenever $\gamma v' = v$, by
$$(z, v,
u).(v', \gamma,t) = (\gamma^{-1}z, v', \phi(z,\gamma)^{-1}tu).$$ (The fact that
this defines an action follows from the definition of $\sigma$ in terms of
$\phi$.) We may now check that these actions commute, since, if $z=\beta y$ and
$\gamma v = v'$,
\begin{align*} [(x,y,v,t)(z,v',u)](v'',\gamma,s) & = (\beta x,\beta v,
\tau(\beta x,z)tu)(v'', \gamma,s) \\ & = (\gamma^{-1}\beta
x,\gamma^{-1}\beta v,
\phi(\beta x, \gamma)^{-1}\tau(\beta x,z)tus), \end{align*} while one has
\begin{align*} (x,y,v,t)[(z, v', u)(v'', \gamma,s)] & = (x,y,v, t)
(\gamma^{-1}z,\gamma v'', \phi(z,\gamma)^{-1}us)\\ & = (\gamma^{-1}\beta x,
\gamma^{-1}\beta v, \tau(\gamma^{-1}\beta
x,\gamma^{-1}z)\phi(z,\gamma)^{-1}tus). \end{align*} The equality of these two
follows from a calculation similar to \cite{CHMM} Lemma 3.2.
\end{proof}

\begin{rems*} 1. Using the orientation reversing diffeomorphism of the Riemann
surface $\Sigma = \Gamma\backslash X$, one can show as in Proposition 7
\cite{CHMM} that the  algebra $C(\Omega) \rtimes_{\overline{\sigma}} \Gamma$ is
isomorphic to
$C(\Omega) \rtimes_{\sigma}\Gamma$, where $\bar\sigma$ denotes the complex
conjugate of $\sigma$.

\noindent 2. Theorem 1.6 shows that the algebras of observables in both the
continuous and discrete models, have the same module structure ($K$-theory),
which is a critical property used to show the equivalence of these two models,
cf. \cite{CHMM}. \end{rems*}

\section{The Discrete Model}

In this section we formulate a version of the integer quantum Hall effect on a
graph in hyperbolic space. The discussion uses a construction due to Sunada
\cite{Sun}, see also \cite{CHMM},
\cite{MM}, \cite{Co2}.

We have seen that the magnetic field defines a multiplier or group 2-cocycle
$\sigma$ on
$\Gamma$. The unperturbed discrete Hamiltonian $H_\sigma$ is then the random
walk operator on the Cayley graph of $\Gamma$, in the projective $(\Gamma,
\sigma)$-representation on $\ell^2(\Gamma)$. More precisely, recall the left
$\sigma$-regular representation
$$ (U(\gamma) f)(\gamma') = f(\gamma^{-1}\gamma') \sigma(\gamma',
\gamma^{-1}\gamma')
$$
$\forall f\in \ell^2(\Gamma)$ and $\forall \gamma, \gamma' \in \Gamma$. It has
the property that
$$ U(\gamma) U(\gamma') = \sigma(\gamma, \gamma') U(\gamma\gamma')
$$ Let $S= \{A_j, B_j, A_j^{-1}, B_j^{-1} : j= 1, \ldots, g\}$ be a symmetric
set of generators for $\Gamma$. Then the unperturbed discrete Hamiltonian is
explicitly given as
$$ H_\sigma :\ell^2(\Gamma) \to \ell^2(\Gamma), \quad H_\sigma =
\sum_{\gamma\in
S} U(\gamma). $$

\begin{defn*} A {\em random} or {\em almost periodic} discrete operator is a
continuous family of bounded operators on $\ell^2(\Gamma)$, $\Omega \ni r
\mapsto T_r \in B(\ell^2(\Gamma))$ which satisfies the following equivariance
condition $$ U(\gamma) T_r U(\gamma)^* = T_{\gamma r} \qquad \forall \gamma\in
\Gamma, \forall r\in \Omega, $$ where the disorder space $\Omega$  is
defined in
section 1.2.
\end{defn*}

Clearly a random (or almost periodic) discrete operator $T$ can be viewed as an
element of the twisted cross product algebra $C(\Omega) \rtimes_\sigma \Gamma$.
The set of all such operators $T$ gives our  class of Hamiltonians in the
discrete model as perturbations: $$ H_{\sigma, T} = H_\sigma + T
$$ This significantly generalises the class of Hamiltonians that were
considered
in \cite{CHMM}. Following
\cite{JiS} we introduce the algebra ${\mathcal R}(\Omega \times \Gamma,
\sigma)$
of functions on $\Omega \times \Gamma$ which are rapidly decreasing in the
$\Gamma$ variable uniformly on $\Omega$.

\begin{lemma}
$f(H_{\sigma, T} ) \in C(\Omega) \rtimes_\sigma \Gamma$ for any continuous
function on
$\R$. If $E\not\in \text{spec}(H_{\sigma, T})$, then $P_E \in {\mathcal
R}(\Omega \times \Gamma, \sigma)$, where $P_E = \chi_{[0,E]}(H_{\sigma, T})$ is
the spectral projection of the random Hamiltonian to energy levels less than or
equal to $E$.
\end{lemma}

\begin{proof} Since $E\not\in \text{spec}(H_{\sigma, T})$, then $P_E =
\chi_{[0,E]}(H_{\sigma, T}) = \phi(H_{\sigma, T})$ for some smooth, compactly
supported function $\phi$. Now by definition, $H_{\sigma, T} \in C_c(\Omega
\times \Gamma, \sigma) \subset {\mathcal R}(\Omega \times \Gamma, \sigma)$, and
since ${\mathcal R}(\Omega \times \Gamma, \sigma)$ is closed under the smooth
functional calculus by the result of \cite{JiS}, it follows that $P_E \in
{\mathcal R}(\Omega \times \Gamma, \sigma)$. \end{proof}

\section{Fredholm module and Connes-Kubo formula: the discrete model} We
restrict our attention in this paper to the discrete model although following
\cite{CHMM} and Section 1, it is possible (but more difficult) to give a
parallel discussion for the continuous model.

\subsection{Canonical cyclic cocycles for the discrete model} We have observed
in Section 2 that the random operator $H_{\sigma, T}$ is in the twisted
algebraic group algebra $C_c(\Omega \times \Gamma, \sigma)$. A spectral
projection into a gap in the spectrum of $H_{\sigma, T}$ is given by the smooth
functional calculus applied to $H_{\sigma, T}$. It follows from \cite{JiS} that
such spectral projections lie in ${\mathcal R}(\Omega \times \Gamma, \sigma)$.
Connes constructs a Fredholm module for
${\mathbb C}\Gamma$ which, following \cite{CHMM} we now adapt to the case of
${\mathcal R}(\Omega \times\Gamma, \sigma)$. A Fredholm module over this
algebra
is given by the graded Hilbert space,
$ L^2(\Omega \times\Gamma) \oplus L^2(\Omega \times\Gamma).$ The grading is the
obvious one given by the $2\times 2$ matrix $\varepsilon
=\left(\begin{array}{cc}
1 & 0 \\ 0 & -1 \end{array}\right)
$ and a (constant) family of operators $F$ parametrised by $\Omega$, which is
taken to be multiplication by the matrix function $\left(\begin{array}{cc} 0 &
\varphi^{*} \\ \varphi & 0 \end{array}\right), $ where we restrict $\varphi$ to
the orbit $\Gamma.u$ in $X$. To describe $\varphi$ first identify
$\ell^2(\Gamma)$ with the $\ell^2$ sections of the restriction of the spinor
bundle to the orbit $\Gamma.u$ in $X$. {}From this point of view
$\varphi$ corresponds to Clifford multiplication of a unit tangent vector to a
geodesic connecting a given vertex of the graph to a point $x_0\notin\Gamma.u$.
We use the same notation $\varphi(\gamma.u)$ for this unit tangent vector,
regarding $\varphi$ as a function from $\Gamma.u$ to $T(\hyp )$, the tangent
space of $\hyp$.

Recall that the {\em area cocycle} $c$ of the group $\Gamma$ is defined as
follows. The group
$PSL(2, \R)$ acts on $\mathbb H$ such that $\mathbb H \cong PSL(2, \R)/SO(2)$.
The area 2-cocycle on $PSL(2, \R)$, cf. \cite{Co2}, is defined by:
$$c(\gamma_1,
\gamma_2) = \text{Area}(\Delta(o, \gamma_1.o, {\gamma_2}^{-1}.o)) \in \R,$$
where $o$ denotes an origin in $\mathbb H$ and
$\text{Area}(\Delta(a,b,c))$ denotes the hyperbolic area of the geodesic
triangle in $\mathbb H$ with vertices at $a, b, c \in \mathbb H$. Then the
restriction of $c$ to the subgroup $\Gamma$ is the area cocycle $c$ of
$\Gamma$.

In \cite{CHMM}, we show using an argument similar to Connes \cite{Co2} that if
$\lambda$ denotes the left regular $\sigma$ representation of
$C^*(\Gamma,\sigma)$ then $[F,\lambda(\gamma)]$ is Hilbert-Schmidt. So
$({\mathcal H}, F)$ is also a 2-summable module for ${\mathbb C}
(\Gamma,\sigma)$. We also determined explicitly the character of this Fredholm
module to be: $$\tr_{c}(f^0,f^1,f^2)= \sum_{\gamma_0\gamma_1\gamma_2=1} f^0(
\gamma_0)f^1( \gamma_1) f^2(
\gamma_2)c(1,\gamma_1,\gamma_1\gamma_2)\sigma(\gamma_1, \gamma_2),$$ for
$f^0,f^1,f^2\in
\C(\Gamma, \sigma)$ and where $c$ is the area 2-cocycle on the group $\Gamma$.
It is known that this cocycle is bounded cf. \cite{Mos}, \cite{Gr} and so by a
well known argument due to \cite{Co2}, the formula above extends to give a
non-trivial element of the cyclic cohomology of the smooth subalgebra
${\mathcal
R}(\Gamma, \sigma)$. In particular, for projections
$P \in {\mathcal R}(\Gamma, \sigma)$, one has \cite{CHMM} $$ \tr_c(P,P,P) =
\Index(PFP) \in 2(g-1) \Z. $$ We shall use this fact in proving the integrality
of other cocycles of interest to us in this paper.

Given a group 2-cocycle $c$ on $\Gamma$, there is a canonical way to define a
cyclic 2-cocycle $\tr_{\Lambda, c}$ on the twisted cross product algebra
$C_c(\Omega \times \Gamma, \sigma)$. First regard the elements of $C_c(\Omega
\times \Gamma, \sigma)$ as $C_c(\Gamma)$ valued functions on $\Omega$. Then for
$f^0,f^1,f^2\in C_c(\Omega \times \Gamma, \sigma)$ $$\tr_{\Lambda,
c}(f^0,f^1,f^2)=
\int_\Omega d\Lambda(r)\tr_{c}(f^0(r),f^1(r),f^2(r))$$ This is merely a
parametrised version of the prescription given in \cite{Co2}, \cite{CHMM} (see
the discussion in Appendix C of \cite{Mc} in the Euclidean case).

\begin{prop} Let $P$ be a projection into a gap in the spectrum of the discrete
random Hamiltonian $H_{\sigma, T}$. Then $P\in {\mathcal R}(\Omega \times
\Gamma, \sigma)$, and may be regarded as a twisted convolution operator by a
(continuously) parametrised function $P$ on $\Gamma$. Then for any $r\in
\Omega$, one has $$
\int_\Omega d\Lambda(r') \tr_c\left(P(r'), P(r'), P(r')\right) = \Index(P(r) F
P(r))
\in 2(g-1)\Z
$$ That is, for any $r\in \Omega$, one has
$$
\langle [\tr_{\Lambda, c}], [P]\rangle = \tr_{\Lambda, c}(P, P, P) =
\Index(P(r)FP(r)) \in 2(g-1)\Z
$$
\end{prop}

\begin{proof} We observe that the LHS of the above equation can be rewritten as
$$
\int_\Omega d\Lambda(r) \sum_{\gamma_0\gamma_1\gamma_2=1}
P(r,\gamma_0)P(r,\gamma_1)P(r, \gamma_2)
c(1,\gamma_1,\gamma_1\gamma_2)\sigma(\gamma_1,\gamma_2) = \int_\Omega
d\Lambda(r) \Index(P(r) F P(r) ) $$ where $P(r)$ denotes the projection $P(r,
\cdot)$, by the result in \cite{CHMM} refered to above. We will now examine the
integer valued function $n(r) = \Index(P(r)FP(r))$ on $\Omega$. Since
$P({\gamma
r}) = U(\gamma) P(r)U(\gamma)^*$, we see that the operators $P({\gamma r}) F
P({\gamma r})$ and $P(r)U(\gamma)^* FU(\gamma) P(r)$ are unitarily equivalent.
In particular
$$
\Index(P({\gamma r}) F P({\gamma r})) =
\Index(P(r)U(\gamma)^* F U(\gamma) P(r)). $$ It is also straightforward to see
that $P(r)FP(r)$ and $P(r)U(\gamma)^* F U(\gamma) P(r)$ differ by a compact
operator and so
$$
\Index(P(r)FP(r)) =
\Index(P(r)U(\gamma)^*F U(\gamma) P(r)). $$ Therefore $$ n(r) = \Index(P(r)F
P(r)) = \Index(P({\gamma r}) F P({\gamma r})) = n(\gamma r)
$$ for all $\gamma \in \Gamma$ and $r\in \Omega$. Note that since the function
$r \to P(r)$ is continuous, it follows that $r \to n(r)$ is continuous. Since
the action of $\Gamma$ on $\Omega$ has a dense orbit, and since $n(r)$ is
continuous, integer valued and $\Gamma$-invariant, it follows that $n(r)$ is a
constant independent of $r\in \Omega$. Finally, by a main result in
\cite{CHMM},
one has $\Index(P(r)F P(r)) \in 2(g-1) \Z$. \end{proof}

\begin{rems} By \cite{JiS}, one knows that ${\mathcal R}(\Omega \times \Gamma,
\sigma)$ is a smooth subalgebra of the twisted cross product algebra
$C(\Omega)\rtimes_\sigma \Gamma$, hence it follows from the proposition above
that the range of the cyclic 2-cocycle $\tr_{\Lambda, c}$ on $K$-theory is: $$
\tr_{\Lambda, c}(K_0({\mathcal R}(\Omega \times \Gamma, \sigma) = \tr_{\Lambda,
c}(K_0(C(\Omega)\rtimes_\sigma \Gamma)) = 2(g-1)\Z. $$ This result has also
been
established as a special case of a topological index formula in \cite{Ma}. See
\cite{Co2} for details on the pairing of cyclic cocycles and $K$-theory.
\end{rems}

\subsection{The hyperbolic Connes-Kubo formula for the conductance: discrete
model}

Our discussion can be viewed as a parametrised version of the approach for
orbifolds given in \cite{MM}. Here we derive the discrete analogue of the
hyperbolic Connes-Kubo formula for the Hall conductance 2-cocycle, in the case
of random Hamiltonians, generalising a result in \cite{MM}. We then relate
it to
the area 2-cocycle on the twisted cross product algebra by the group $\Gamma$
given in the previous subsection, and we show that these define the same cyclic
cohomology class. This enables us to use the results of the previous section to
show that the Hall conductance is integral, having plateaus at all energy
levels
belonging to any gap in the spectrum of the Hamiltonian.

The Cayley graph of the group $\Gamma$, which acts freely on hyperbolic space
$X$ embeds in the hyperbolic plane as follows. Fix a base point $u \in X$ and
consider the orbit of the $\Gamma$ action through $u$. This gives the vertices
of the graph. The edges of the graph are geodesics constructed as follows. Each
element of the group $\Gamma$ may be written as a word of minimal length in the
generators of $\Gamma$ and their inverses. Each generator and its inverse
determines a unique geodesic emanating from a vertex $x$ and these form the
edges of the graph. Thus each word $x$ in the generators determines a piecewise
geodesic path from $u$ to $x$.

We will now describe the hyperbolic Connes-Kubo formula for random
Hamiltonians.
Let $\Omega_j$ denote the (diagonal) operator on $\ell^2(\Gamma)$ defined by
$$
\Omega_j f(\gamma) = \Omega_j(\gamma)f(\gamma) \quad \forall f\in
\ell^2(\Gamma)
\quad \forall \gamma\in \Gamma
$$ where
$$
\Omega_j(\gamma) = \int_o^{\gamma.o} a_j \quad j=1,\ldots ,2g $$ and where $\,
\{a_j\} \quad j=1,\ldots ,2g$ is the lift to $\mathbb H$ of a symplectic basis
of harmonic 1-forms on the Riemann surface of genus $g$ associated to $\Gamma$.

For $j=1,\ldots ,2g$, define the derivations $\delta_j$ on ${\mathcal R}(\Omega
\times \Gamma, \sigma)$ as being the commutators $\delta_j a = [\Omega_j,
a]$. A
simple calculation shows that
$$
\delta_j a (r, \gamma) = \Omega_j(\gamma) a(r, \gamma) \quad \forall a\in
{\mathcal R}(\Omega \times \Gamma, \sigma) \quad \forall \gamma\in \Gamma. $$
Note that these are not inner derivations, and also that we have the simple
estimate
$$ |\Omega_j(\gamma)|\le ||a_j||_{(\infty)} d(\gamma.o, o) $$ where
$d(\gamma.o,
o)$ and the distance in the word metric on the group $\Gamma$, $d_\Gamma
(\gamma, 1)$ are equivalent. This then yields the estimate
$$ |\delta_j a (r, \gamma)| \le C_N d_\Gamma (\gamma, 1)^{-N}\quad \forall N\in
\mathbb N
$$ i.e $\delta_j a \in {\mathcal R}(\Omega \times \Gamma, \sigma) \quad \forall
a \in {\mathcal R}(\Omega \times \Gamma, \sigma)$. Note that since $\forall
\gamma,
\gamma' \in \Gamma$, the difference
$\Omega_j(\gamma\gamma') - \Omega_j(\gamma')$ is a constant independent of
$\gamma'$, we see that $\Gamma$-equivariance is preserved. For $j=1,\ldots
,2g$,
define the cyclic 2-cocycles $$
\tr^K_j (a_0,a_1,a_2) =
\int_\Omega d\Lambda (r) \tr (a_0(r) (\delta_ja_1(r) \delta_{j+g} a_2(r) -
\delta_{j+g} a_1 (r) \delta_j a_2 (r) )) $$ where $a_0, a_1, a_2 \in {\mathcal
R}(\Omega \times \Gamma, \sigma)$ These have been argued in section 6,
\cite{CHMM} to give the Hall conductance for currents in the $(j+g)$th
direction
which are induced by electric fields in the $j$th direction. Then the {\em
hyperbolic Connes-Kubo formula} for the Hall conductance for random potentials
is the cyclic 2-cocycle given by the sum $$ \tr^K_\Lambda (a_0,a_1,a_2) =
\sum_{j=1}^g \tr^K_j (a_0,a_1,a_2).
$$

\begin{thm}[The Comparison Theorem]
$$ [\tr^K_\Lambda] = [\tr_{\Lambda, c}] \in HC^2({\mathcal R}(\Omega \times
\Gamma, \sigma) $$
\end{thm}

\noindent{\bf Sketch of Proof:} Our aim is now to compare the two cyclic
2-cocycles and to prove that they differ by a coboundary i.e. $$
\tr^K_\Lambda (a_0,a_1,a_2) - \tr_{\Lambda, c} (a_0,a_1,a_2) = b\lambda_\Lambda
(a_0,a_1,a_2)
$$ for some cyclic 1-cochain $\lambda_\Lambda$ and where $b$ is the cyclic
coboundary operator. This is easily deduced from results in \cite{MM} where we
observe that $\lambda_\Lambda (a_0,a_1) = \int_\Omega d\Lambda (r) \lambda
(a_0(r), a_1(r))$.

It follows from Connes' pairing theory of cyclic cohomology and $K$-theory
\cite{Co2} and the Comparison Theorem above that

\begin{cor} For all projections
$P\in {\mathcal R}(\Omega \times \Gamma, \sigma)$, one has $$
\tr^K_\Lambda (P, P, P) = \tr_{\Lambda, c} (P, P, P) \in 2(g-1)\Z. $$
\end{cor}

Finally, we think of our discrete model as describing electrons hopping along
the vertices of the Cayley graph of $\Gamma$ via the geodesics joining them. In
a very strong magnetic field which is uniform and normal in direction to the
sample and at low temperatures quantum mechanics dominates so that our
Hamiltonian models the dynamics of these electrons in the presence of random
impurities in the Cayley graph. (One should picture the graph embedded in the
hyperboloid.)

The hyperbolic Connes-Kubo formula for the Hall conductance $\sigma_E$ at the
energy level $E$ is defined as follows; let $P_E = \chi_{[0,E]}(H_{\sigma, T})$
be the spectral projections of the Hamiltonian to energy levels less than or
equal to $E$. Then if $E\not\in \text{spec}(H_{\sigma, T})$, one can show that
$P_E \in {\mathcal R}(\Omega \times \Gamma, \sigma)$, and the Hall conductance
is defined as
$$
\sigma_E = \tr^K_\Lambda (P_E, P_E, P_E). $$ When this is combined with the
previous corollary, one sees that the Hall conductance takes on values in
$2(g-1) \mathbb Z$ whenever the energy level $E$ lies in a gap in the spectrum
of the random Hamiltonian
$H_{\sigma, T}$. In fact we notice that the Hall conductance is a integer
valued
constant function of the energy level $E$ for all values of $E$ in the same gap
in the spectrum of the random Hamiltonian. We now summarize this discussion:

\begin{thm}[Hyperbolic quantum Hall effect in the presence of disorder $\#1$]
Suppose that the energy level $E$ lies in a gap in the spectrum of the random
Hamiltonian $H_{\sigma, T}$, then the Hall conductance $$
\sigma_E = \tr^K_\Lambda (P_E, P_E, P_E) = \tr_{ \Lambda, c} (P_E, P_E,
P_E) \in
2(g-1) \mathbb Z
$$ That is, the Hall conductance has plateaus which are {\em integer} multiples
of $2(g-1)$ on any gap in the spectrum of the random Hamiltonian. \end{thm}

\begin{rems*} This theorem should be compared with Theorem 4.6 in the next
section, where we relax the continuity assumption on the projection at the
expense of imposing the more stringent condition that the $\Gamma$-action is
ergodic on the disorder space $\Omega$.
\end{rems*}

\section{Gaps in extended states and the main theorem}

In this Section, we extend our results in Section 3 to consider gaps in
extended
states by imposing the condition that the $\Gamma$ action on $\Omega$ admits a
quasi-$\Gamma$-invariant ergodic measure. We firstly show that the domain of
definition of the conductance cocycle is in fact an algebra which is larger
than
the algebra of observables that was discussed earlier. This enables us to study
the conductance for random Hamiltonians which do not necessarily have any gaps
in their spectrum. We also show that there are only a finite number of gaps in
extended states for any random Hamiltonian which corresponds to what is
expected
to be physically observed. In contrast, it is sometimes possible to have an
infinite number of gaps in the spectrum of a random Hamiltonian, cf.
\cite{CHMM},
\cite{CEY},
\cite{Rief}.

\subsection{Preliminaries} We begin by recalling some standard notions on
Fredholm modules
$({\mathcal{H}},F)$ for a $*$-algebra $A$.

\noindent{\bf Definition 1} {\em Hilbert-Schmidt algebra extensions}

Let $A$ be a $*$-algebra and $({\mathcal{H}}, F)$ be a 2-summable Fredholm
module over $A$ i.e.
\begin{enumerate}
\item there is a $*$-representation $\pi$ of $A$ in $\mathcal{H}$. \item
$F=F^*,
F^2=1$ and $[F, \pi(a)]\in {\mathcal{L}}^{2} \quad\forall\ a\in A$
\end{enumerate} where ${\mathcal{L}}^{2}$ denotes the $*$-ideal of Hilbert
Schmidt operators. Define the $*$-subspace of {\em Hilbert-Schmidt}
elements  of
the von Neumann closure
$W_A$ of $A$ as
\[ W^2_A = \{f\in W_A\ | \quad [F, \pi(f)]\in {\mathcal{L}}^{2}\} \] where we
observe that the $*$-representation $\pi$ extends to the von Neumann closure
$W_A$. Since
${\mathcal{L}}^{2}$ is a $*$-ideal, and since $[F, \cdot]$ is a derivation, one
sees that $W^2_A$ is a $*$-algebra. It then follows that $A \subset W^2_A
\subset W_A$ are $*$-subalgebras and that
$({\mathcal H},F)$ is also a 2-summable Fredholm module over $W^2_A$ by
definition. In particular, its Chern character $\Ch({\mathcal H},F)$ also
extends to $W^2_A$ and one has
$$ {\mathbb Z} \ni \Index(\pi(P)F\pi(P)) = \left\langle\Ch({\mathcal H},F),
[P]\right\rangle =
\tr \left(\pi(P) [F,\pi(P)] [F,\pi(P)]\right)\quad\forall\ P\in\text{ Proj
}(W^2_A)
$$ Recall that the Besov space $B^{1/2}_2$ is the space of bounded measurable
functions on $\R$ such that
\[
\int\int |f(x+t) - f(x)|^2 t^{-2} dxdt<\infty \] Then it follows from a theorem
of Peller that $W^2_A$ is closed under the Besov space $B^{1/2}_2$ functional
calculus, cf. \cite{Co2}.

\noindent{\bf Definition 2} {\em Sobolev algebra extensions}

In the notation above, define the $*$-algebra of Sobolev elements in $W_A$
as \[
W^{2,\delta}_A = \{f\in W_A | \quad \pi(\delta_j(f)) \in {\mathcal L}^2\quad
\forall\ j=1,\dots,2g\}
\] Since $\delta_j$ is a derivation and ${\mathcal L}^2$ is a $*$-ideal, it
follows that $W^{2,\delta}_A$ is a $*$-algebra. It is also clear that $A
\subset
W^{2,\delta}_A \subset W_A$ are $*$-subalgebras.

Note that the intersection $W^2_A \cap W^{2,\delta}_A$ is also a $*$-subalgebra
of $W_A$ which contains $A$. When $A= {\mathcal R}(\Omega\times\Gamma, \sigma)$
(we will assume this for the rest of the section) and $(\mathcal H, F)$ is the
Fredholm module considered in the previous section, we observe that the cyclic
2-cocycles $\tr_{\Lambda,c}$ and $\tr^K_\Lambda$ are well defined on the
intersection $W^2_A \cap W^{2,\delta}_A$. Moreover on $W^2_A \cap
W^{2,\delta}_A$ one still has the equality $\tr_{\Lambda,c}-\tr^K_\Lambda
=b\lambda_\Lambda$. In particular one sees that
$$
\tr_{\Lambda, c}(\pi(P), \pi(P),\pi(P)) = \tr^K_\Lambda (\pi(P), \pi(P),\pi(P))
= \Index(\pi(P)F\pi(P)) \in 2(g-1){\mathbb Z} $$ $\forall\ P\in\text{ Proj
}(W^2_A \cap W^{2,\delta}_A)$. Notice that the algebra $W^2_A \cap
W^{2,\delta}_A$ is much larger than $A$, being closed under the Besov space
functional calculus and not merely the continuous functional calculus, thus we
have considerably extended the results that we obtained in the previous
Section.
We formulate this in terms of a theorem in the next subsection.

\subsection{Gaps in extended states}

The gaps in extended states of the random Hamiltonian $H_{\sigma, T}$
correspond
precisely to the plateaus of the Hall conductance in the quantum Hall effect.
One can make this precise as follows. Let $P_E = \chi_{(-\infty, E]}(
H_{\sigma,
T}) $ denote a spectral projection and $\sigma(a_0, a_1, a_2) =
\tr_\Lambda(\epsilon a_0 da_1 da_2)$ is the conductance cyclic 2-cocycle, where
$da = [F, a]$ and $\epsilon$ denotes the grading operator. We also use the
notations $\sigma_P = \sigma(P,P,P)$ and $\sigma_E = \sigma_{P_E}$. We define
the {\em domain} of the conductance 2-cocycle $\sigma_E$, $\rm{domain}(\sigma)$
to be the set of all $E \in \R$ such that $P_E$ is both a Sobolev projection as
well as a Hilbert-Schmidt projection in $W_A$ i.e the set of all $E\in \R$ such
that
$P_E \in W^2_A \cap W^{2,\delta}_A$. The picture that motivates our definition
of extended states is the following:

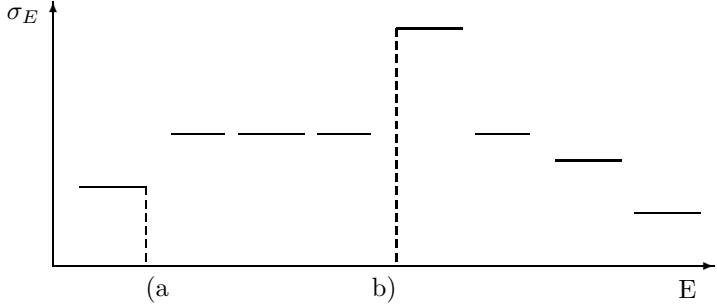
\begin{picture}(300,120)(-50,-10)
\put(0,0){\vector(1,0){250}}
\put(0,0){\vector(0,1){100}}
\put(10,30){\line(1,0){25}}
\put(45,50){\line(1,0){20}}
\put(70,50){\line(1,0){25}}
\put(100,50){\line(1,0){20}}
\put(130,90){\line(1,0){25}}
\put(160,50){\line(1,0){20}}
\put(190,40){\line(1,0){25}}
\put(220,20){\line(1,0){25}}
\multiput(35,2)(0,5){6}{\line(0,1){3}}
\multiput(130,2)(0,5){18}{\line(0,1){3}} \put(35,-4){\makebox(0,0)[tl]{(a}}
\put(130,-4){\makebox(0,0)[tr]{b)}}
%\put(38,-3){\makebox(0,0)[t]{a}}
%\put(127,-3){\makebox(0,0)[t]{b}}
\put(240,-5){\makebox(0,0)[t]{E}}
\put(-5,95){\makebox(0,0)[r]{$\sigma_E$}} \end{picture}

Hence we define a {\em gap in extended states} of the Hamiltonian $H_{\sigma,
T}$ to be an open interval $(a, b)$ of $\R$ which satisfies the following
conditions. \begin{enumerate}
\item $(a, b)\cap \rm{domain}(\sigma) \ni E \mapsto \sigma_E$ is a constant
function on the set $(a, b)\cap \rm{domain}(\sigma)$; \item $(a, b)$ is maximal
in the sense that if $(a', b')$ is another open interval in $\R$ satisfying
condition 1 above, and which has nonempty intersection with $(a, b)$, then
$(a',
b')
\subset (a, b)$. \end{enumerate}

Then $E \in
\R$ belongs to a gap in extended states of the Hamiltonian $H_{\sigma, T}$ if
the spectral projection
$P_E = \chi_{(-\infty, E]}( H_{\sigma, T}) $ is both a Sobolev projection as
well as a Hilbert-Schmidt projection in $W_A$. The motivation for this
definition comes from the fact that a gap in the spectrum of the Hamiltonian
$H_{\sigma, T}$ can be also described as an interval $(a, b)$ such that
\begin{enumerate}
\item $(a, b)\ni E \to N_E = \tr_\Lambda(P_E)$ is a constant function on the
interval $(a, b)$;
\item $(a, b)$ is maximal in the sense that if $(a', b')$ is another interval
satisfying condition 1 above, and which has nonempty intersection with $(a,
b)$,
then $(a', b') \subset (a, b)$. \end{enumerate}

Our first result in this Section is to prove that there are only a {\em finite}
number of gaps in extended states of the random Hamiltonian $H_{\sigma, T}$. To
achieve this we introduce a couple of auxilliary cyclic 2-cocycles, which are
defined as \begin{align*} \psi(a_0, a_1, a_2) & = -\int_\Omega
d\Lambda(r)\tr(a_0(r) da_1(r) da_2(r)) \\ \xi^\pm(a_0, a_1, a_2) & =
-\int_\Omega d\Lambda(r)\tr((1\pm \epsilon) a_0(r) da_1(r) da_2(r))
\end{align*} where $da = [F, a]$ and $\epsilon$ denotes the grading
operator. We
also recall the notion of positivity for cocycles, \cite{Co2}. A cyclic
2-cocycle $\eta$ is said to be {\em positive} if the formula $$ \eta(b_0^*a_0,
a_1, b_1^*) = \langle a_0 \otimes a_1, b_0\otimes b_1 \rangle $$ defines a
positive sesquilinear form on the vector space $W_A\otimes W_A$.

We begin with the following lemma.
\begin{lemma} The cyclic 2-cocycles
$\psi(a_0, a_1, a_2), \xi^\pm(a_0, a_1, a_2) $ are positive \end{lemma}
\begin{proof} We just have to show that
$$
\psi(a_0^*a_0, a_1, a_1^*) \ge 0 \mbox{ and } \xi^\pm(a_0^*a_0, a_1, a_1^*) \ge
0.
$$ We verify this for $\psi$ first. Observe that since $F^* = F$, one has
$da_1^* = [F, a_1^*] = - [F, a_1]^* = -(da_1)^*$. Then \begin{align*}
\psi(a_0^*a_0, a_1, a_1^*) & = -\int_\Omega d\Lambda(r)\tr(a_0(r)^*a_0(r)
da_1(r) da_1(r)^*)\\ & = -\int_\Omega d\Lambda(r) \tr(a_0(r) da_1(r)
da_1(r)^*a_0(r)^*)\\ & = \int_\Omega d\Lambda(r) \tr(a_0(r) da_1(r)
(da_1(r))^*a_0(r)^*)\\ & = \int_\Omega d\Lambda(r) \tr(a_0(r) da_1(r)(a_0(r)
da_1(r))^*) \ge 0. \end{align*} One can verify that $\xi^\pm$ is positive in a
similar manner to the above, by observing that $2(1\pm \epsilon) = (1\pm
\epsilon)^2$.
\end{proof}

We recall \cite{Co2} that any cyclic 2-cocycle $\eta$ is {\em additive} in the
sense that if $P_1$ and $P_2$ are projections in $W^2_A \cap W^{2,\delta}_A$
which are orthogonal to each other, then one has
$$\eta(P_1+P_2,P_1+P_2,P_1+P_2)
= \eta(P_1,P_1,P_1) + \eta(P_2,P_2,P_2).$$ We use the notations $\psi_P =
\psi(P,P,P), \quad \xi^\pm_P = \xi^\pm(P,P,P)$ and $\psi_E = \psi_{P_E}, \quad
\xi^\pm_E = \xi^\pm_{P_E}$. Then as an immediate consequence of the positivity
and additivity properties of the cyclic 2-cocycles $\psi$ and $\xi^\pm$ one has
\begin{cor} The functions
$E \to \psi_E$ and $E\to \xi^\pm_E $ are monotone increasing on their
domains of
definition.
\end{cor}

\begin{cor} The following inequality holds $$\psi_P \ge |\sigma_P|$$ for all
projections $P$ in the intersection of domains of definition of the two
functions. \end{cor}

\begin{proof} Just observe that we have the identities $\psi - \sigma = \xi^+$
and $\psi + \sigma = \xi^-$.
\end{proof}

\begin{thm}[Gaps in extended states theorem] There are only a {finite}
number of
gaps in extended states of the random Hamiltonian $H_{\sigma, T}$ with disorder
space having an ergodic
$\Gamma$ action.\end{thm}
\begin{proof} Let $E_1 < E_2$ be energy levels that are in two different
gaps in
extended states of the random Hamiltonian $H_{\sigma, T}$. Then since the
spectral projections are related by $P_{E_2} = P_{E_1} + P$, where the
projection $P= P_{E_2} - P_{E_1} \in W^2_A \cap W^{2,\delta}_A$, it follows by
hypothesis that $\sigma_P \ne 0$. By additivity of $\psi$ and corollary
4.3, one
has
$$\psi_{E_2} =
\psi_{E_1} + \psi_P \ge \psi_{E_1} + |\sigma_P|.$$ That is, by Proposition 4.5
below, $\psi_{E_2}
\ge \psi_{E_1} +2(g-1) $. By induction, we conclude that if $E_1<E_2 \ldots
E_{n-1}<E_n$ are energy levels that lie in distinct gaps in extended states of
the random Hamiltonian $H_{\sigma, T}$, then one has $$\psi_{E_{j+1}} \ge
\psi_{E_j} + 2(g-1)$$ for all $j=1, \ldots, n-1$. Also by corollary 4.2, the
function $E \to \psi_E$ is a positive, increasing function on its domain of
definition. It follows that $\psi_{E_n} \ge \psi_{E_n} - \psi_{E_1} \ge
2(g-1)(n-1)$. So the number of gaps in extended states is less than or equal to
the integer part of $\frac{1}{2(g-1)}\psi_E +1$, where $E\le ||H_{\sigma, T}||$
is the largest real number that lies in a gap in extended states of the random
Hamiltonian $H_{\sigma, T}$.
\end{proof}

We now prove an analogue of Proposition 3.1 for gaps in extended states.
Observe
first that
$W_A=L^\infty(\Omega)\times_\sigma\Gamma$ so that {\em a priori} any
$P\in\text{
Proj }(W^2_A \cap W^{2,\delta}_A)$ is given by a measurable operator valued
function $r\to P(r)$ (acting on $\ell^2(\Gamma)$). Hence the function $n(r)$
introduced in the proof of proposition 3.1 is only measurable and so we need an
additional assumption. The one we choose is to assume that the $\Gamma$-action
on $\Omega$ admits a quasi-invariant ergodic measure. This implies that $r\to
n(r)$ is $\Gamma$ invariant and measurable and hence constant $\Lambda$ a.e.
Therefore we have

\begin{prop} Let $P$ be a projection into a gap in extended states of the
discrete random Hamiltonian $H_{\sigma, T}$ with disorder space having an
ergodic
$\Gamma$ action. Then $P\in \text{ Proj }(W^2_A \cap W^{2,\delta}_A)$, and may
be regarded as a twisted convolution operator by a (measurably) parametrised
function $P$ on $\Gamma$. Then for almost any $r\in\Omega$, one has
$$
\int_\Omega d\Lambda(r') \tr_c\left(P(r'), P(r'), P(r')\right) = \Index(P(r) F
P(r))\in 2(g-1)\Z.
$$ That is, for almost any $r\in \Omega$, one has $$
\langle [\tr_{\Lambda, c}], [P]\rangle = \tr_{\Lambda, c}(P, P, P) =
\Index(P(r)FP(r)) \in 2(g-1)\Z.
$$
\end{prop}

We can now state the main result of the paper, which follows from Proposition
4.5 and from the discussion in section 4.1.

\begin{thm}[Hyperbolic quantum Hall effect in the presence of disorder:$\#2$]
Suppose that the energy level $E$ lies in a gap in extended states of the
random
Hamiltonian $H_{\sigma, T}$ with disorder space having an ergodic
$\Gamma$ action, then the Hall conductance $$
\sigma_E = \tr^K_\Lambda (P_E, P_E, P_E) = \tr_{ \Lambda, c} (P_E, P_E,
P_E) \in
2(g-1) \mathbb Z.
$$ That is, the Hall conductance has plateaus which are {\em integer} multiples
of $2(g-1)$ on any gap in extended states of the random Hamiltonian. \end{thm}

\end{document}